\documentclass[12pt]{notices}
\usepackage{graphicx}
\usepackage{amssymb}
\usepackage{epstopdf}
\usepackage{amsmath,amscd}
\usepackage{amsthm}
\usepackage{url,verbatim}
\usepackage{stmaryrd}
\usepackage{caption,subcaption}

\usepackage[dvipsnames]{xcolor}

\newcommand\EE{{\mathbb E}}

\newcommand\RR{{\mathbb R}}
\newcommand\ZZ{{\mathbb Z}}

\newcommand\PP{{\mathbb P}}

\newcommand\TT{{\mathbb T}}

\newcommand\pc{p_{\text{\rm c}}}
\newcommand\pslab{p_{\text{\rm slab}}}
\renewcommand\o{\text{\rm o}}
\newcommand\eqd{\mathrel{\stackrel{\text{\rm d}}{=}}}

\newcommand\oo{\infty}

\newcounter{mycount1}\newcounter{mycount2}\newcounter{mycount3}\newcounter{mycount}

\numberwithin{equation}{section}
\numberwithin{theorem}{section}
\numberwithin{figure}{section}
\newcommand\sym{\,\triangle\,}

\title{Harry Kesten (1931--2019)\\ A personal and scientific tribute}

\author{
Geoffrey R. Grimmett
  \affil{Geoffrey Grimmett is Professor Emeritus of Mathematical Statistics at
   Cambridge University, and Professor of Applied Probability at the University of Melbourne. His email address is  grg@statslab.cam.ac.uk.}
  \and
Gregory F. Lawler
  \affil{Gregory Lawler is the George Wells Beadle Distinguished Service Professor at the
  University of Chicago. His email address is lawler@math.uchicago.edu.}
}

\begin{document}

\date{\today}

\maketitle

\begin{figure}
\centerline{\includegraphics[width=0.45\textwidth]{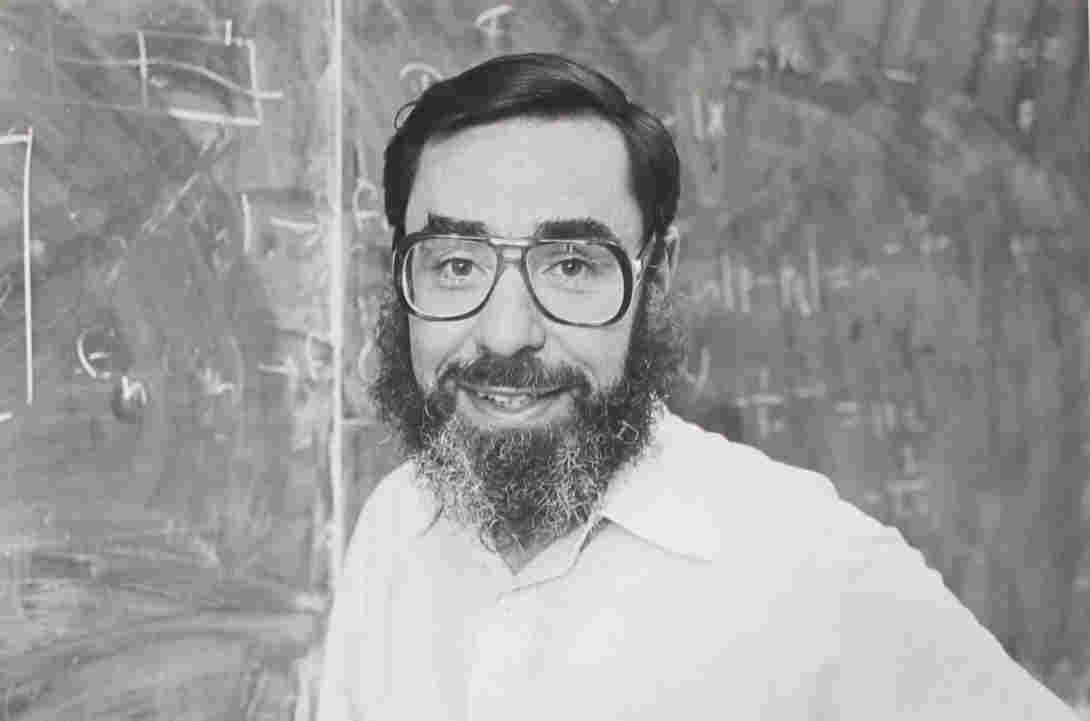}}
\caption*{{\bf Figure 1.} Harry Kesten in White Hall, Cornell University, around
1980.}
\end{figure}

The mathematical achievements of Harry Kesten since the mid-1950s have revolutionized 
probability theory  as a subject in its own right and in
its associations with aspects of algebra, analysis, geometry, and statistical physics.
Through his personality and scientific ability, he has framed the modern subject 
to a degree exceeded by no other.

Harry inspired high standards of honesty, modesty, and informality, and 
he played a central part in the creation of a lively and open community of researchers.  

\section*{Biographical notes}

\subsection*{Early life}

Harry Kesten's early life was far from tranquil.
His already migrant father migrated once again from Germany following Hitler's appointment as Chancellor in 1933.
Harry survived the war under the protection of a non-Jewish family in the Netherlands. 
His parents died naturally during and immediately following the war.   
In 1952, he recognised his likely future as a mathematician.

Harry was born in Duisburg, Germany on 19 November  1931 to
Michael and Elise Kesten. \lq\lq Hamborn (Hochfeld)" is listed as his place of birth
in both a German and  a  Polish document.  Harry had Polish citizenship  from birth
via his father, and he retained this until his American naturalization in 1962.

He was the only child of Michael Kesten 
(b.\ Podo{\l}owoczyska, \cite{podw}, then in Galicia and now Ukraine, 
22 November 1890) and Elise Abrahamovich (b.\ Charlottenburg, Berlin, 4 October 1905).
The Kesten family name featured prominently in the affairs of the 
substantial Jewish population of  Podo{\l}owoczyska
in the late 19th century. Michael moved from Galicia to Berlin, where he and Elise were married in
Charlottenburg on 23 May 1928.
Amos Elon has written eloquently in \cite{elon} of the Jewish community in Berlin up to 1933,
the year in which the Kesten family moved to Amsterdam,
perhaps as a member of a group of Jewish families. 

German military forces attacked and invaded the Netherlands in May 1940.
Shortly afterwards, Harry was offered protection by
a non-Jewish Dutch family residing in Driebergen near Utrecht, with whom he lived until the Dutch liberation
in May 1945.  His father was hidden in the same village, and they could be
in occasional contact throughout the occupation. During that period Harry attended school
in a normal way (\cite{dynkin}). 
Elise died of leukaemia in Amsterdam
on 11 February 1941, and Michael died in October 1945, probably in Groningen of  cancer.
Later in life, Harry kept in touch with his Dutch family, and would visit them whenever he was in the Netherlands.
  
Harry moved back to Amsterdam at the end of the war to live with an older married cousin
who had been born a Kesten, and who had survived the war in Switzerland. 
It was during the period between 1945 and graduation from
High School in 1949 that his attachment to Orthodox Judaism developed, and this was to remain
with him for most of the rest of his life. 

The earliest surviving indication of Harry's scientific ability is found in his school report 
on graduation from the  General Secondary High School in 1949. His six grades in the six given scientific topics
(including mathematics)
are recorded as  
five 10s and one 9. In languages, they were
7 (Dutch), 7 (French), 8 (English), 10 (German). There were three subjects in which his mark was a
mere 6 (\lq satisfactory'), including handwriting and physical education. In later life, Harry was very active physically,
and was keen to run, swim, hike, and to ski cross-country, usually with friends and colleagues. 

\begin{figure}[h]
\centerline{\includegraphics[width=0.45\textwidth]{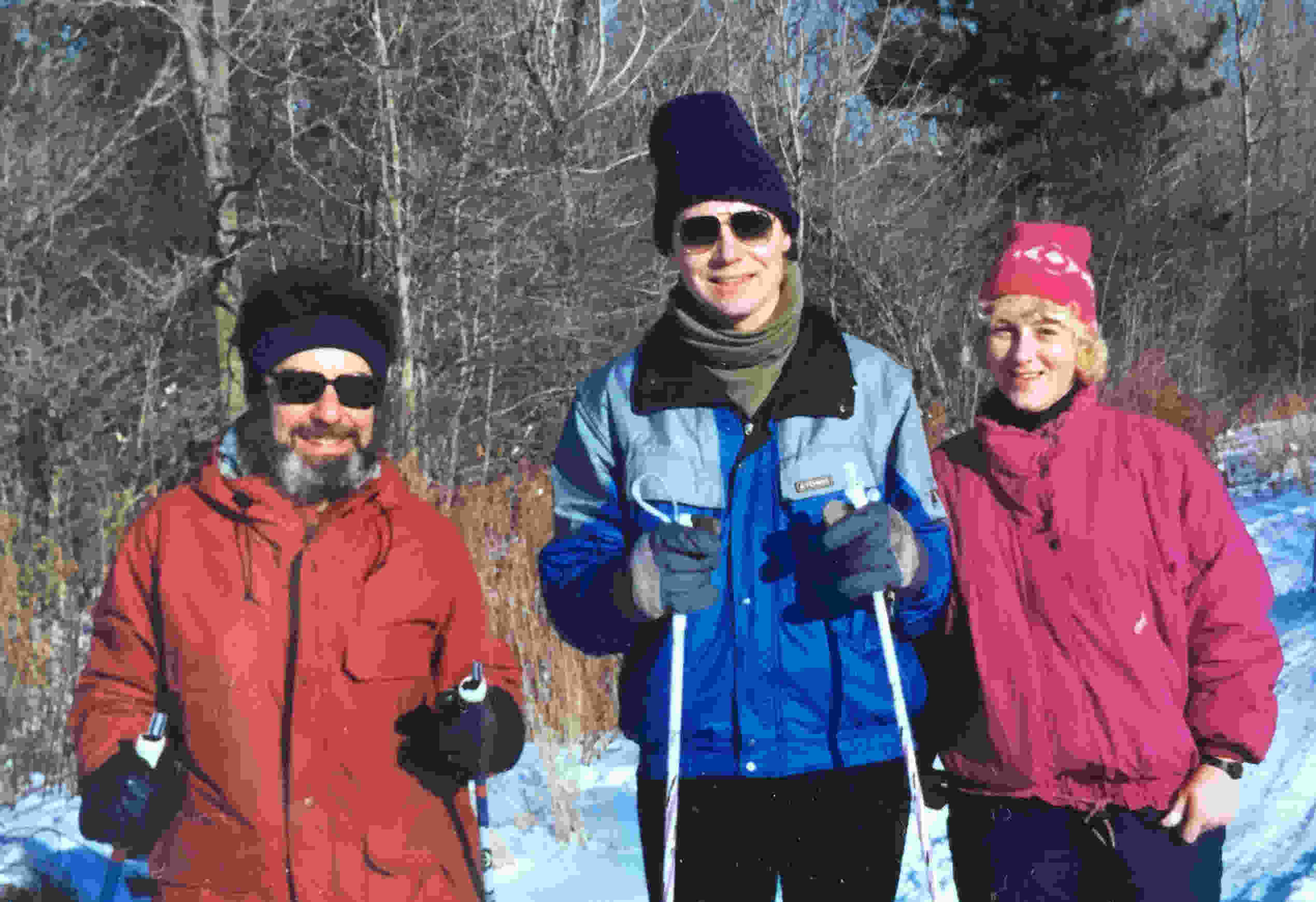}}
\caption*{{\bf Figure 2.} 
Harry Kesten, cross-country skiing
with Rob and Margriet van den Berg, Ithaca, 1991.}
\end{figure}

Following his uncle's advice to become a chemical engineer, 
Harry entered the University of Amsterdam in 1949 to study chemistry (\cite{dynkin}).
This was not altogether successful, and Harry took a  particular dislike to laboratory work. 
He moved briefly to theoretical
physics before settling on mathematics.
From 1952 to 1956 he had a half-time assistantship in
the Statistical Department of the Mathematical Centre (now the CWI), Amsterdam,
under the supervision of David van Dantzig (known for his theory of collective marks)
and Jan Hemelrijk.
He shared an office with fellow student Theo (J.\ Th.) Runnenburg,
with whom he wrote his first papers on topics in renewal theory and queueing theory.
The pair of papers \cite{HK-first} are notable, since
they are probably related to the Master's (almost, in a sense, doctoral) thesis that Harry wrote in 1956.

It was around this time that he met his wife-to-be, Doraline 
Wabeke, who worked in the Mathematical Centre Library while studying interior design at evening school.

\subsection*{Middle years} 
Mark Kac visited Holland in 1955, and Harry had the opportunity to meet him at the Mathematical Centre. 
He wrote to Kac in January 1956 
to enquire of a graduate fellowship at Cornell University to study probability theory, perhaps for one year.  
Van Dantzig wrote to Kac in support \lq\lq \dots I have not the slightest doubt  that,
if you grant him a fellowship, you will consider the money well spent afterwards.''

A Junior Graduate Fellowship was duly arranged with a stipend of \$1400 plus fees,  
and Harry joined the mathematics graduate program at Cornell that summer,
traveling on a passport issued by the International Refugee Organization.
His fellowship was extended to the next academic year 1957--58 with support from Mark Kac: 
\lq\lq Mr.\ Kesten is $\dots$ the best student we have had here in the last twenty years.
$\dots$ one of these days we will indeed be proud of having helped to educate an
outstanding mathematician.''

He defended his PhD thesis at the end of that year, on the (then) highly novel topic
of random walks on groups. This area of Harry's creation remains an active 
and fruitful area of research at the time
of this memoir.

Doraline followed Harry to the USA in 1957 under the auspices of \lq The Experiment of International Living',
and took a position in Oswego, NY, about 75 miles north of Ithaca on Lake Ontario. As a result of the Atlantic crossing 
she became averse to long boat journeys, and never traveled thus
again.  She studied and converted to Judaism with a rabbi in Syracuse and the couple was married in 1958.

Harry and Doraline moved to Princeton in 1958, where  Harry held
a (one-year) instructorship in the company of Hillel Furstenberg. 
For the following academic year, he accepted a position at the Hebrew
University of Jerusalem. Harry was interested in settling in Israel, and wanted to try it out, 
but there were competing pressures from Cornell, who wished to attract him back to Ithaca, and from Doraline's
concerns about practical matters. 
They postponed a decision on the offer
of an Assistant Professorship at Cornell (with a standard nine hours/week teaching load), opting instead
to return for the year 1961--62 on a one-year basis.

It was during that year that Harry and Doraline decided to make Ithaca their home.
Harry was promoted to the rank
of Associate Professor in 1962, and in 1965 Harry and Doraline celebrated both his promotion to  Full Professor 
and the birth of their only child, Michael. Harry stayed at Cornell for his entire career, becoming Emeritus Professor in 2002.

Harry, Doraline, and Michael lived for many years with their cats at 35 Turkey Hill Road,
where visitors would be welcomed for parties and walks.

Unsurprisingly, many invitations arrived, and he would invariably try to oblige.  
Harry traveled widely, and paid many extended visits to universities and research centers in the USA and abroad, 
frequently accompanied by Doraline, and in earlier years by Michael.

For some years the principal events in Cornell probability included the biweekly 5 mile runs
of Harry with Frank Spitzer, and able-bodied visitors were always welcome. When Harry's knees 
showed their age, he spent time swimming lengths in Cornell's Teagle pool,
usually in the now \emph{pass\'e} men-only sessions.  He worked on his problems while swimming.

\begin{figure}[h]
\centerline{\includegraphics[width=0.3\textwidth]{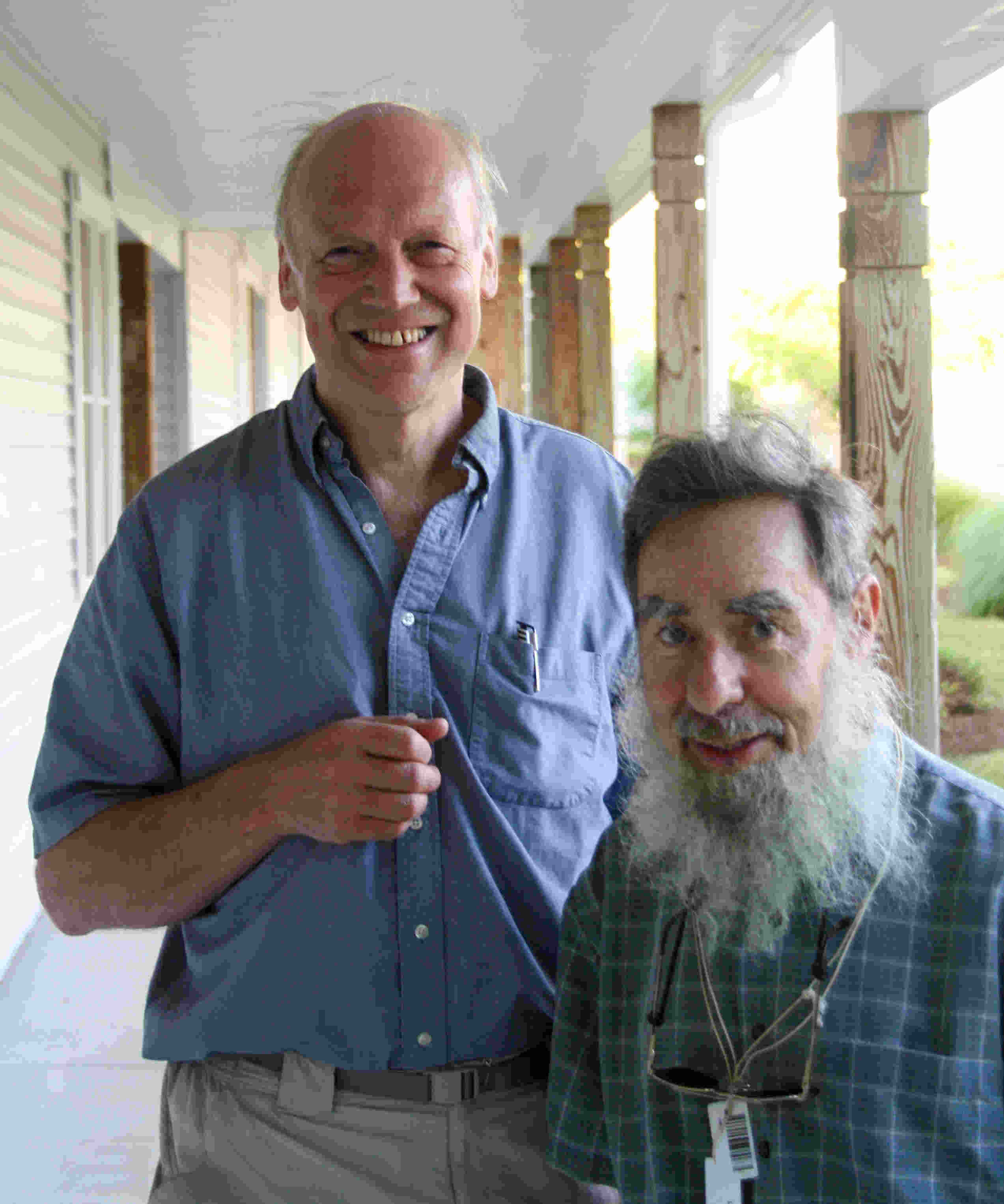}}
\caption*{{\bf Figure 3.}
Harry Kesten with Geoffrey Grimmett in Kendal at Ithaca, 2011.}
\end{figure}
\subsection*{Later years}

Harry maintained his research activities and collaborations beyond his retirement from active duty in 2002.
He spoke in the Beijing ICM that year on the subject of percolation, finishing with a slide listing individuals
who had been imprisoned in China for the crime of expressing dissent (\cite{icm2002}).

He was awarded an honorary doctorate at the Universit\'e de Paris-Sud in 2007,
shortly following his diagnosis with Parkinson's disease.
Harry and Doraline sold their house and moved in 2008 into the Kendal retirement home,
a \lq home from home' for numerous retired Cornell academics. 
Doraline developed Alzheimer's disease and died at Kendal on 2 March 2016,
followed by Harry from complications of Parkinson's disease on 29 March 2019.
A volume on percolation remained on his bedside table until the end.

\subsection*{Awards of distinction}

From among the awards made to Harry Kesten, mention is made here of
the Alfred P.\ Sloan Fellowship (1963), the Guggenheim Fellowship (1972), the Brouwer Medal
(1981), the SIAM George P\'olya Prize (1994), and the AMS Leroy P.\  Steele Prize for Lifetime
Achievement (2001).
Harry delivered the Wald Memorial Lectures of the IMS (1986), and was elected a Correspondent
of the Royal Netherlands Academy of Arts and Sciences (1980),
a Member of the National Academy of Sciences (1983), and of
the American Academy of Arts and Sciences (1999).
He was elected an Overseas Fellow of Churchill College, Cambridge (1993), and
was awarded an Honorary Doctorate of the Universit\'e de Paris-Sud (2007).

He was an invited lecturer at three ICMs, Nice (1970), Warsaw (1983), and Beijing (2002), and he spoke
at the Hyderabad ICM (2012) on the work of Fields Medalist Stanislav Smirnov. He was a
member of the inaugural class of Fellows of the AMS in 2013.

\section*{Personality and influence}

Probability theory gained great momentum
in the second half of the 20th century. Exciting and beautiful problems were formulated and solved,
and connections with other fields of mathematics and science, both physical and socio-economic, were established.
The general area attracted a large number of distinguished scientists, and it grew
in maturity and visibility. 
Harry was at the epicenter of the mathematical aspects of this development. He contributed 
new and often startling results at the leading edge
of almost every branch of probability theory. 

Despite an occasionally serious aspect,
he was a very sociable person who enjoyed his many scientific collaborations and was a popular
correspondent. His archive of papers (now held by Cornell University)
reveals a wealth of letters exchanged with many individuals worldwide, and
every serious letter received a serious reply, frequently proposing solutions to the problems posed. 
He was especially keen to discuss and collaborate with younger people, and he played a key role
through his achievements  and personality in bringing them into the field.

Harry commanded enormous respect and affection amongst those who knew him well. He 
displayed an uncompromising honesty, tempered by humanity,  in 
both personal and professional matters.
This was never clearer than in his opposition to oppression,
and in his public support for individuals deprived of their
positions, or even liberty, for expressing their beliefs or needs. 

Harry loved hard problems. Supported by an extraordinary technical ability and a total lack of fear,  
he gained a just reputation as a fearsome problem-solver. His work often exceeded the greatest current expectations,
and it could be years before the community caught up with him. When the going became too
tough for the rest of us, he would simply refuse to give in. The outcome is 196 works listed on MathSciNet,
almost every one of which contains some new idea of substance.  
An excellent sense of Harry as a mathematician may be
gained by reading the first two pages of Rick Durrett's appreciation 
\cite{MR1703122}, published in 1999 to mark $40+$ years of Harry's mathematics.

\section*{Scientific work}

Most areas of probability theory have been steered, even moulded, by Harry, and
it is not uncommon to attend conferences in which a majority of speakers refer to his 
work as fundamental to their particular topics.
In this memoir, we do not aim at a comprehensive survey but, instead, to
select and describe some high points.  Our selection is personal by necessity,
and readers desirous of a  more comprehensive account are referred to \cite{MR1703122,grimmett}. 

\begin{figure}
\centerline{\includegraphics[width=0.225\textwidth]{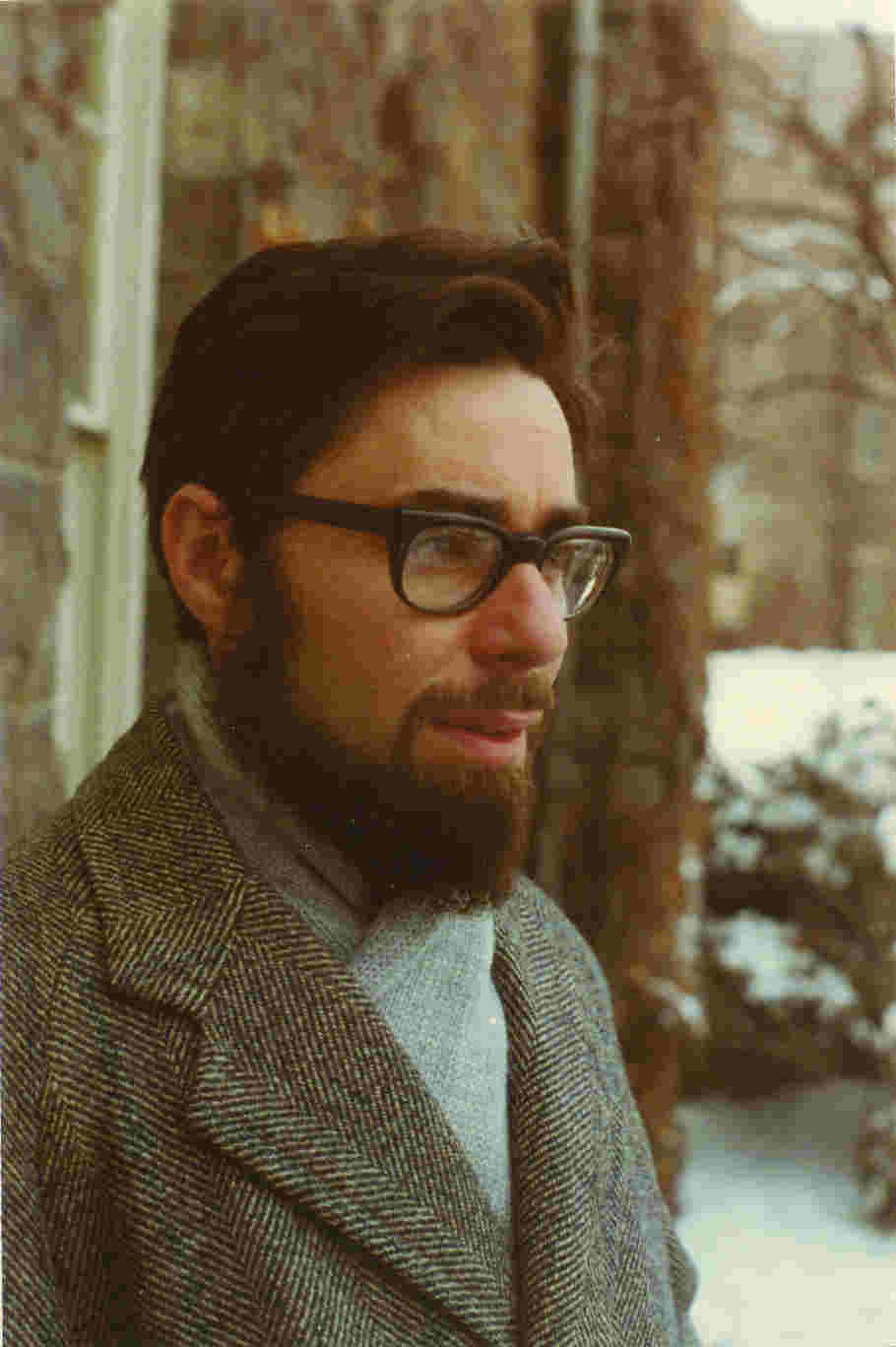}\quad
\includegraphics[width=0.22\textwidth]{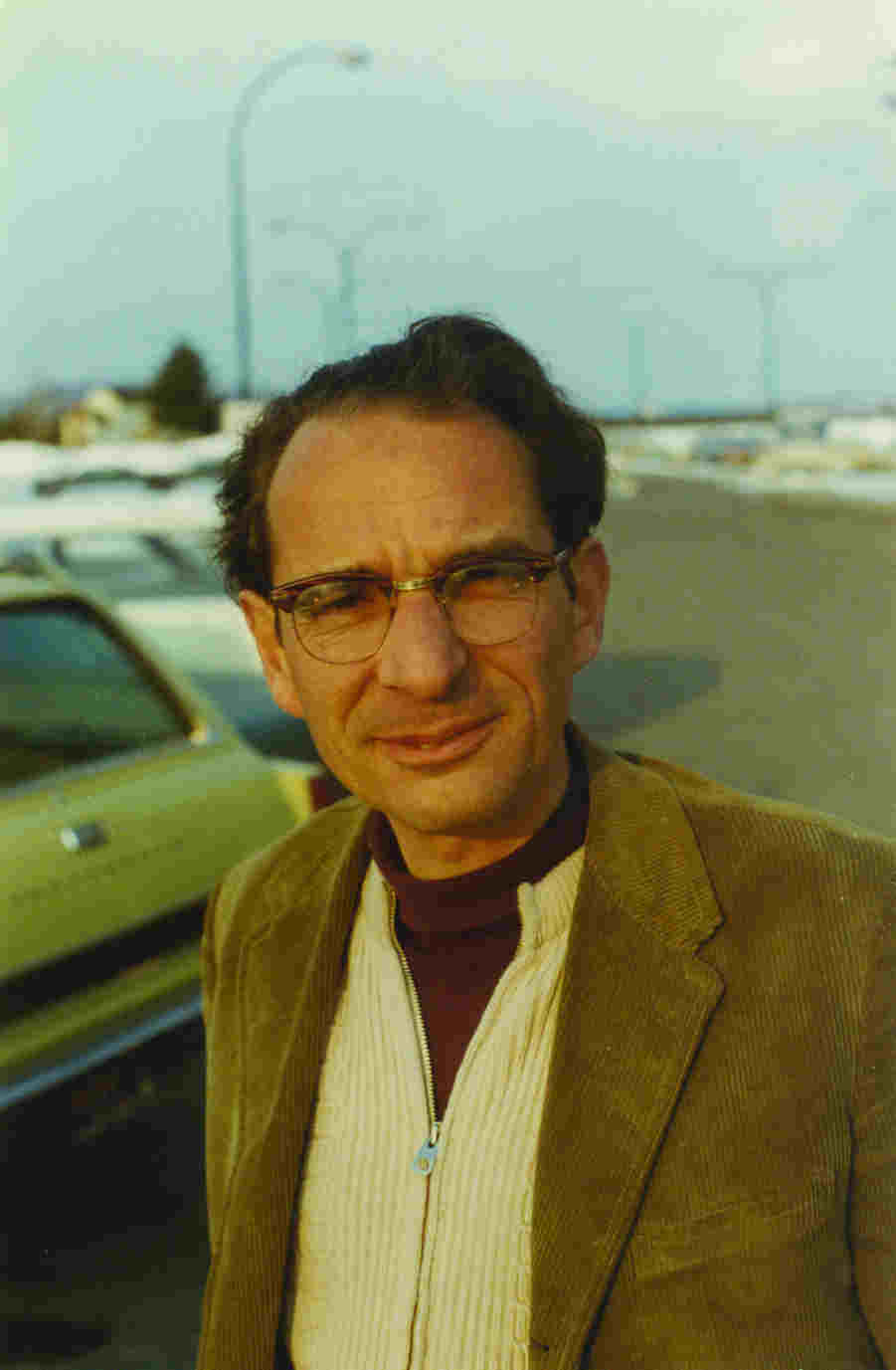}}
\caption*{{\bf Figure 4.}
Harry Kesten and Frank Spitzer in 1970.}
\end{figure}

\subsection*{Random walks on groups}

Harry's PhD work was on the topic of random walks
on groups.  So-called  \lq simple random walk' takes random jumps of size $\pm 1$ about the
line $\ZZ$. The domain $\ZZ$ may replaced by either $\ZZ^d$ or $\RR^d$,
and the unit moves replaced by a general family of independent and identically distributed displacements.
Harry's space was algebraic rather than Euclidean.
He considered a countable group $G$, a symmetric ($p(x) = p(x^{-1})$) probability distribution
on a generating set of $G$, and defined random walk on $G$
as the process that moves at each step from $y$ to $yx$
with probability $p(x)$.  Let $q_{2n}$ denote the
probability that the random walk returns to its starting point after
$2n$ steps.  For example, the usual
random walk on $\ZZ^d$ has $q_{2n} \sim c_d  n^{-d/2}$, while
for the free group on two generators we have $q_{2n}\approx\bigl(\frac34\bigr)^n$. 

Harry showed for a general countably infinite group $G$ that the quantity
\[       \Lambda(G,p):= \lim_{n \rightarrow \infty} q_{2n}^{1/2n} \]
equals both the spectral radius and the maximal value
of the spectrum of the associated operator
on $L^2(G)$ given by the random walk. He proved that the equality
$\Lambda(G,p)=1$ is a property of the group $G$
and not of the particulars of the transition  probabilities $p$, and if it holds we write
$\Lambda(G)=1$.  The so-called
\lq Kesten criterion for amenability' states that 
$\Lambda(G) =1$ if and only if $G$ is amenable.\footnote{Of the various 
equivalent definitions of amenability, the reader is reminded that
a discrete group $G$ is amenable if there exists a sequence $F_n$ of finite subsets such that, for $g\in G$,
$|(g F_n)\sym F_n|/|F_n| \to 0$ as $n\to\oo$.}
This remarkable characterization of amenability may be viewed as a fairly early contribution
to the currently important area of geometric group theory (see \cite{MR3329732}).

\subsection*{Products of random matrices}  

One of the earliest papers in the now important field of \emph{random matrices} 
is \cite{MR121828} by Furstenberg and Kesten, 
written by two
Princeton instructors in 1958--59. They
were motivated by the 1954 work of Bellman, studying the asymptotic behavior of the product
of $n$ independent, random $2 \times2$ matrices, and they derived substantial extensions of Bellman's results.
This now classical paper \cite{MR121828} has been very influential and is much cited, 
despite having proved unwelcome at 
the authors' first choice of journal. It deals with products of random matrices, in contrast to most of
the modern theory which is directed towards spectral properties.

Consider a stationary sequence $X_1,X_2,\dots$ of random $k \times k$ matrices, and let 
$$
Y^n=(y^n_{ij})=X_nX_{n-1}\cdots X_1.
$$
In an analysis termed by Bellman \lq\lq difficult and ingenious'',
Furstenberg and Kesten proved a law of large numbers and a central limit theorem. Firstly,
if $X$ is ergodic, the limit $E:= \lim_{n\to\oo} n^{-1}\log \|Y^n\|_1$ exists a.s.
Secondly, subject to certain conditions,
the limit of $n^{-1}\EE(\log y^n_{ij})$ exists, and  
$n^{-1/2}\bigl(\log y^n_{ij} - \EE(\log y^n_{ij})\bigr)$ is asymptotically normally distributed.

Although they used subadditivity in the proofs, they did not anticipate
the forthcoming theory of subadditive stochastic processes, initiated in 1965 by Hammersley and Welsh 
to study first-passage percolation, 
which would one day provide a neat proof of some of their results.

Harry returned in 1973 to a study of products of random matrices arising in
stochastic recurrence relations. In the one-dimensional case,
it was important to understand the tail behavior of
a random variable $Y$ satisfying a stochastic equation of the form  $Y\eqd MY+Q$;
that is, for random variables $M$ and $Q$
with given distributions, $Y$ and $MY+Q$ have the same distribution.
He proved in \cite{MR440724} that such $Y$ are generally heavy-tailed in that $\PP(Y>y)$
decays as a power law as $y\to\oo$.  This work has generated a 
very considerable amount of interest since in probability, statistics, and mathematical finance.

\subsection*{Random walks and L\'evy processes}

Harry's early years at Cornell marked a heyday for the theory of random walk.  
His colleague Frank Spitzer had written 
his Springer monograph \textit{Principles of Random Walk}, and, together, they and 
their collaborators did much to further the field.  The classical  definition of random walk is as a sum
$ S_n = x + X_1 + \dots + X_n $
where $x$ is an initial position and  $X_1,X_2,\dots$ are independent and identically
distributed.  When the distributions are nice, the theory parallels  
that of the continuous potential theory of the Laplacian.  
When the hypotheses on distributions are weakened, some but not all such properties persist.  

In an early sequence of papers, Harry proved a family of ratio limit theorems for 
probabilities associated with random walks on $\ZZ^d$. 
Here are two examples, taken from joint work with
Spitzer and Ornstein. 
Let $T$ denote the first return time of $X$ to its starting point.
Then $\PP_0(T>n+1)/\PP_0(T>n)\to 1$ as $n\to\oo$, and this may be used
to show that the limit 
$$  
a(x) = \lim_{n\to\oo} \frac{\PP_x(T > n)}{\PP_0(T > n)}
$$ 
exists and equals the potential kernel (fundamental solution of the corresponding discrete
Laplacian) at $x$, that is
$$
a(x) = \sum_{n=0}^\oo \bigl[\PP_0(X_n=0)-\PP_x(X_n=0)\bigr],\quad x \ne 0.
$$
This result requires no further assumptions on the random walk. 

A L\'evy process is a random process in continuous time with stationary independent increments.
Let $X$ be a L\'evy process on $\RR^d$.
The fundamental question arose through work of Neveu, Chung, Meyer, and McKean of deciding when the hitting
probability $h(r)$ by $X$ of a point  $r\in\RR^d$ satisfies $h(r)>0$. 
Harry solved this problem in his extraordinary AMS Memoir of 1969, \cite{MR0272059}. 
The situation is simplest when $d=1$, for which case Harry showed that (apart from special cases)
$h(r)>0$ if and only if the so-called characteristic exponent of $X$ satisfies a certain integral condition.

The expression `random walk' is sometimes used loosely 
in the context of interesting and challenging 
problems, arising for instance in  physical and 
biological models, which lack the full assumptions of independence and identical distribution
of jumps. 
Harry responded to the challenge to attack many of these difficult
problems for which the current machinery was not sufficient.   
He combined his mastery of classical techniques with 
his `problem-solving' ability to develop new ideas for such novel topics. 

A significant variant of the classical random walk is the `random
walk in random environment' (RWRE).  In the one-dimensional case, this is 
given by  (i) sampling random variables for each site, that prescribe the transition probabilities when one reaches the
site, and (ii) performing a random walk (or more precisely a
Markov chain) with those transition probabilities.  RWRE 
is a Markov chain \emph{given the environment}, but it is not itself 
Markovian because, in observing the process, one accrues information
about the underlying random environment.

One of the first RWRE cases
considered was a nearest neighbor, one-dimensional walk, sometimes called a birth--death chain on $\ZZ$.
Let $\alpha_x$ denote the probability that the random walk
moves one step rightwards when at position $x$ (so that $\beta_x =1-\alpha_x $ is the
probability of moving one step leftwards).
We assume the $\alpha_x$ are independent and identically distributed.
In the deterministic environment with, say,  $\beta_x = \beta <\frac12$ for all $x$, 
then as $n \rightarrow \infty$,  $X_n/n
\rightarrow 1 - 2\beta $; when $\beta = \frac12$, the walk returns to the origin
infinitely often.  For the RWRE (with random $\alpha_x$), 
the properties of the ratio $\beta_0/\alpha_0$ are pivotal for determining whether or not
$X_n \to\oo$; a straightforward
birth--death argument indicates that $X_n \rightarrow \infty$ 
if $\EE[\log (\beta_0/\alpha_0)] < 0$.
The regime with $ \EE[\log (\beta_0/\alpha_0)] < 0$ and $\EE[\beta_0/\alpha_0] > 1$ turns out to be
interesting.  It was already known for this regime that $X_n\to\oo$ but $X_n/n\to 0$.
With Kozlov and
Spitzer \cite{kks}, Harry showed in this case
that the rate of growth of $X_t$ is essentially determined by the value of the parameter $\kappa$ defined by
\[   \EE\left[(\beta_0/\alpha_0)^\kappa\right]= 1.\]
In particular, if $\kappa<1$, then $X_n$ has order 
$n^{\kappa}$.
It turns out that the limit distribution of $n^{-\kappa}X_t$ may be expressed
in terms of a stable law with index $\kappa$. This result
answered a question of Kolmogorov.

The above is an early contribution to the broad area
of RWRE, and to the the related
area of homogenization of differential operators with
random coefficients.
The one-dimensional case is special (roughly speaking, because the walk
cannot avoid the exceptional regions in the environment),  and rather precise results are 
now available, including Sinai's proof that $X_n/(\log n)^2$ converges weakly (to a distribution
calculated later by Harry).   RWRE in higher dimensions poses a very challenging problem.
This process is more diffusive than its one-dimensional cousin, and its study requires
different techniques. Major progress has been made on it by Bricmont, Kupiainen, and 
others.

\subsection*{Self-avoiding walk}
A \emph{self-avoiding walk} (or \lq SAW') on the $d$-dimensional hypercubic lattice $\ZZ^d$ is a
path that visits no vertex more than once. SAWs may be viewed as a simple model for long-chain polymers,
and the SAW problem 
is of importance in physics as in mathematics.

Let $\chi_n$ be the number of $n$-step SAWs starting at the origin.
The principal SAW problems  are to establish the asymptotics of $\chi_n$ as $n\to\oo$, and to determine the typical
radius of a SAW of length $n$. Hammersley and Morton proved
in 1954 that there exists a \lq connective constant' $\kappa$ such that $\chi_n=\kappa^{n(1+\o(1))}$ as $n\to\oo$. 
The finer asymptotics of $\chi_n$ have proved elusive, especially in dimensions $d=2,3$.

Harry's 10 page paper \cite{MR152026} from 1963 contains a theorem and a technique; 
the theorem is essentially unimproved, 
and the technique has frequently been key to the work of others since. His main result is the ratio limit theorem
that $\chi_{n+2}/\chi_n \to \kappa^2$ as $n \to\oo$. It remains an open problem to prove that
$\chi_{n+1}/\chi_n\to\kappa$. His technique is an argument now referred to \lq Kesten's pattern theorem'.
To paraphrase Frank Spitzer from  Mathematical Reviews, the idea is that any configuration of $k$ steps which can 
occur more than once in an $n$-step SAW has to occur at least $an$ times, for some $a>0$, 
in all but \lq very few' such SAWs. 

The pattern theorem is proved using a type of path surgery that has been useful in numerous other contexts since.
The proof is centered around an \lq exponential estimate', of a general type that made powerful
appearances in various different settings in Harry's later work.

One of the most prominent current conjectures in SAW theory is that SAW in two dimensions converges,
in an appropriate limit, to a certain Schramm--Loewner evolution (namely, SLE$_{8/3}$).  If one could show that
the limit exists and exhibits conformal invariance, then it
is known that the limit must be SLE$_{8/3}$.
Although the behavior of this SLE is now 
understood fairly precisely (\cite{LSWsaw}), and the
analogues of the finite asymptotics of $\chi_n$ and the 
typical length of a SAW are known for the continuous model,
the problems of establishing that 
the discrete SAW has a limit, and showing that the limit is
conformally invariant,    remain wide open.
 
\begin{figure}[h]
\centerline{\includegraphics[width=0.4\textwidth]{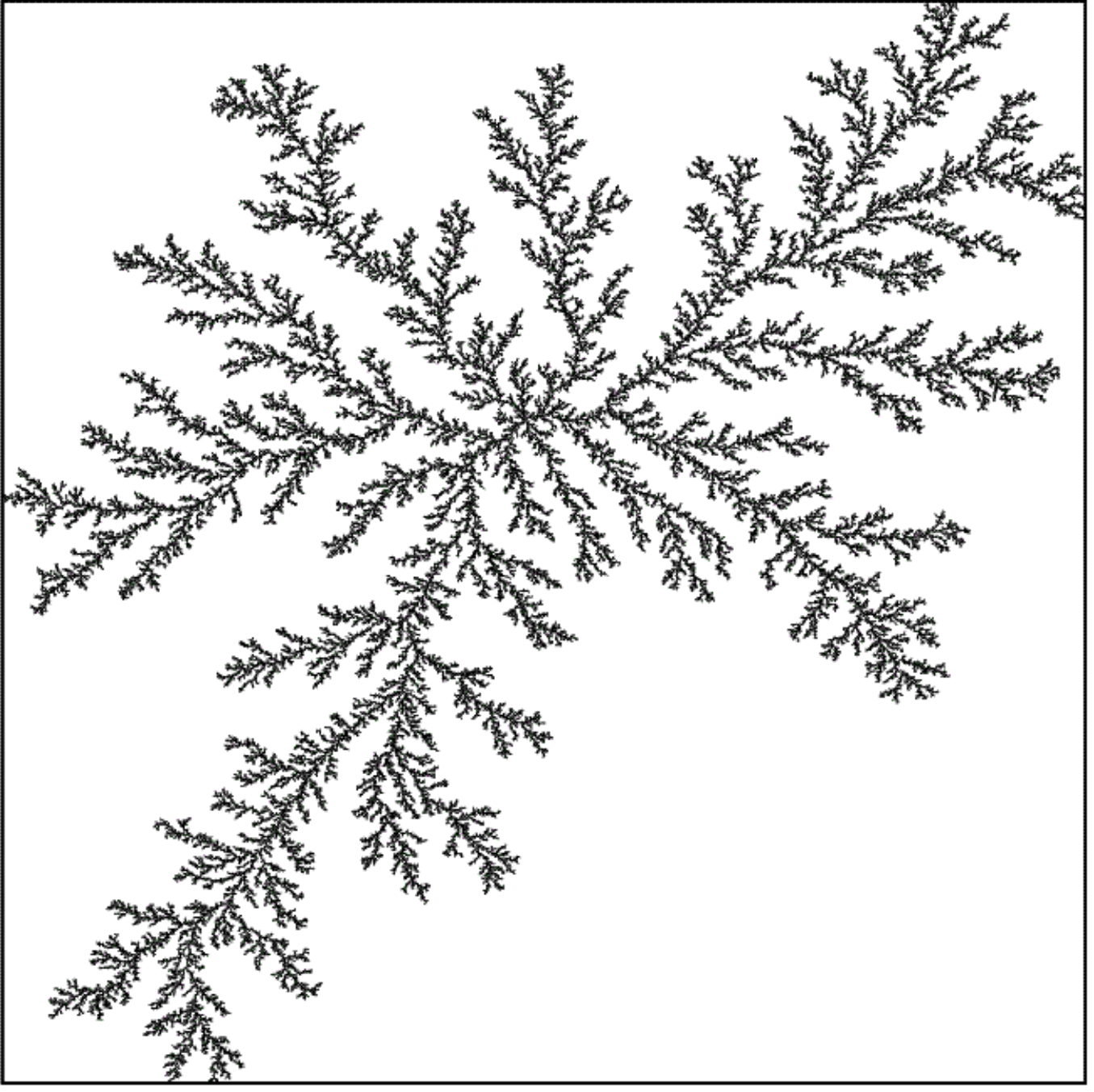}}
\caption*{{\bf Figure 5.} A simulation of diffusion limited aggregation in two dimensions.}
\end{figure}

\subsection*{Diffusion limited aggregation}

Diffusion limited aggregation (DLA) is a growth model introduced by Witten and Sander.   
The model may be defined in general dimensions $d$, but we concentrate here on the case $d=2$.

Let $A_0$ be the origin of the square lattice $\ZZ^2$. Conditional on the set
$A_n$, the set $A_{n+1}$ is obtained by starting a
random walk `at infinity' and stopping it when it reaches
a point that is adjacent to $A_n$, and then adding that
new point to the existing $A_n$.  The concept  of random walk
from infinity can be made precise using harmonic 
measure, and the hard problem is then to describe the evolution of the
random sets $A_n$.  Computer simulations suggest a
random, somewhat tree-like, fractal structure for $A_n$ as $n\to\oo$.  Indeed,
assigning a `fractal dimension' to such a set is a
subtle issue, but a start is made by trying to find the exponent $\alpha$
such that the diameter of $A_n$ grows like $n^{\alpha}$.  Since $A_n$ is a connected set of $n+1$ points, we have
the trivial bounds $\frac12 \leq \alpha \leq 1$.

Harry wrote a short paper \cite{kdla} containing a beautiful
argument showing that $\alpha \leq \frac23$.  For many readers,
this seemed a good start to the problem, and much subsequent effort 
has been invested in seeking improved estimates.  
Unfortunately, no one has yet made
a substantial rigorous improvement to Harry's bound.  There are
a number of planar models with diffusive limited growth, and it is open whether they are in the same universality class. 
On the other hand, \lq Kesten's bound' is one property they have in common.

Harry's argument is simple, but in order to complete it,
he needed a separate lemma (\cite{kbeurling}) about planar random 
walks which is a discrete analogue of a theorem from complex variables due to Beurling.  Kesten's lemma has 
itself proved an extremely useful tool over the last thirty years in the development of the
theory of conformally invariant limits of two-dimensional random walks and other processes.  
Although the original theorem of Beurling was phrased in terms of complex analysis, both  the 
continuous version of Beurling and the discrete lemma of Kesten become most important in
\emph{probabilistic} approaches to problems.  
While Kesten figured out what the right answer should be in his discrete version, 
he was fortunate to have a colleague, Clifford Earle, who knew Beurling's result and 
could refer Harry to a proof that proved to be adaptable (with work\,!) to the 
discrete case.

\subsection*{Branching processes}

The branching process (sometimes called the Galton--Watson process)
is arguably the most fundamental stochastic model for population growth.
Individuals produce a random number of offspring.  
The offspring produce offspring similarly, and so on, with different family-sizes being independent
and identically distributed.
The first question is whether or not the population survives forever.  A basic fact taught in elementary courses on stochastic processes states that the population dies out with probability $1$ if and only if $\mu \leq 1$, where $\mu$ is the mean number of offspring per individual. (There is a trivial deterministic exception to this, 
in which each individual produces exactly one offspring.)

One of Harry's best known results is the Kesten--Stigum theorem \cite{Kstigum} for the supercritical
case $\mu > 1$; it was proved in the more general situation with more than one type of individual, but we
discuss here the situation with only one population type.  If $X_n$
denotes the population size of the $n$th generation, we have $\EE[X_n] = \mu^n X_0$, 
and with some (computable) probability $q >0$, the population survives forever.  
One might expect that, for large $n$, $X_n \sim K_\oo \mu^n$ 
with $K_\oo$ a random 
variable determined by the growth of the early generations.  
In other words, once the population has become large, we should be able to approximate 
its growth by the deterministic
dynamics $X_{n+1}   \sim  \mu X_n$.  If this were true,
we would write $K_\infty = \lim_{n\to\oo} K_n$ where $K_n=  X_n/\mu^n$.
The process $K_n$ is a martingale, and the martingale
convergence theorem implies that $K_n $ converges
almost surely to some limit $K_\oo$.  
However, it turns out to be possible that $K_\oo=0$ a.s., even though the process survives forever
with a strictly positive probability. 

Let $L$ be a random variable with the family-size
distribution, so that $\mu = \EE[L]$. 
It was a problem of some importance to identify a necessary and sufficient moment condition for
the statement that $\EE[K_\infty] = X_0$. 
Kesten and Stigum showed that this condition is that  $\EE[L \log^+L] < \infty$.

Harry was inevitably attracted  by the \emph{critical} branching process (with $\mu=1$).
Amongst his numerous results are the well known
necessary and sufficient conditions (proved with Ney and Spitzer)
for the so-called Yaglom and Kolmogorov laws,
$$
\PP(X_n>0) \sim \frac cn,\quad \PP(X_n>nx\mid X_n>0) \to e^{-c'x}.
$$
We return later to his work on random walk on the critical family 
tree conditioned on non-extinction.

\subsection*{Percolation}

\begin{figure}
\centerline{\includegraphics[width=0.47\textwidth]{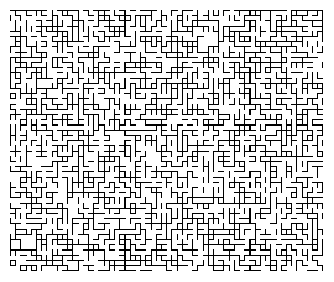}}
\caption*{{\bf Figure 6.}
A simulation of bond percolation on $\ZZ^2$ with $p=0.51$.}
\end{figure}

The percolation model for a disordered medium was pioneered by  Hammersley in the 1950/60s,
and it has since become one of the principal objects in probability theory.
In its simplest form, each edge of the square lattice is declared \emph{open} with probability $p$ and
otherwise \emph{closed}, different edges having independent states. How does the geometry of the open 
graph vary as $p$ increases, and in particular for what $p$ does there exist an infinite open cluster? 

It turns out that there exists a \emph{critical probability} $\pc$
such that an infinite open cluster exists if and only if $p>\pc$. The so-called \lq  phase transition'
at $\pc$ is emblematic of phase transitions in mathematical physics.
Percolation theory is frequently used directly in the study of other systems,
and it has led to the development of a number of powerful insights and
techniques.      
 
\begin{figure}
\centerline{\includegraphics[width=0.45\textwidth]{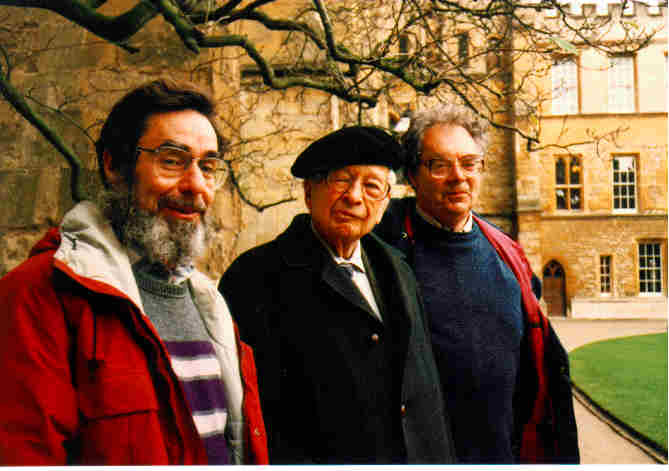}}
\caption*{{\bf Figure 7.}
Harry Kesten, Rudolf Peierls, and Roland Dobrushin in New College, Oxford, 1993.}
\end{figure}
 
As Harry wrote in the preface of his 1982 book \cite{hk82},
\begin{quote}
 \lq\lq Quite apart from the fact that percolation theory had its origin in an honest applied problem $\dots$, 
it is a source of fascinating problems of the best kind a mathematician can
wish for: problems which are easy to state with a minimum of preparation, but whose 
solutions are (apparently) difficult and require new methods.
At the same time many of the problems are of interest to or proposed by statistical physicists 
and not dreamt up merely to demonstrate ingenuity."
\end{quote}
They certainly require ingenuity to solve, as demonstrated in Harry's celebrated proof  that
$\pc=\frac12$ for the square lattice problem, published in 1980 (\cite{MR575895}, see
Figure 6). Harry
proved that $\pc\le \frac12$, thereby complementing Harris's earlier proof that $\pc\ge \frac12$.
His paper, and the book \cite{hk82} that followed, resolved this notorious open problem, and
invigorated an area that many considered almost impossibly mysterious. 

Harry's book \cite{hk82} was a fairly formidable work containing many new results 
for percolation in two dimensions, set in quite a general context. 
He was never frightened by technical difficulty or complication, and entertained
similar standards of the reader.  This project led in a natural way to
his important and far-sighted work \cite{MR879034} on scaling relations 
and so-called \lq arms' at and near the critical point, that was to prove relevant in 
the highly original study initiated by Schramm and developed by Smirnov 
and others on conformal invariance in percolation.

Write $C$ for the open cluster containing the origin.
According to scaling theory, macroscopic functions, such as the percolation probability $\theta$
and mean cluster size $\chi$, given by
$$
\theta(p)=\PP_p(|C|=\oo), \quad \chi(p)=\EE_p|C|,
$$ 
have singularities at $\pc$ of the form $|p-\pc|$ raised to an appropriate power called a \lq critical exponent'.
In similar fashion, when $p=\pc$, several random variables 
associated with the open cluster at the origin have power-law tail behaviors of the form $n^{-\delta}$
as $n \to\oo$, for suitable critical exponents $\delta$. The set of critical exponents describes 
the nature of the singularity and they are characteristics of the model. They are not, however, independent
variables, in that they satisfy the \lq scaling relations' of statistical physics. It is an open problem
to prove almost any of the above in general dimensions. 

In the special case of two dimensions, 
the proof of existence of critical exponents had to wait beyond the invention of SLE by Schramm around 2000,
and the proof of Cardy's formula by Smirnov in 2001 (illustrated in Figure 8), 
and is the work of several individuals including
Lawler and  Werner. In a precursor  \cite{MR879034} of that, Harry
proved amongst other things that, conditional on the existence of certain exponents, certain others
must also exist and a variety of scaling relations ensue.
In this work, he introduced a number of techniques that have been at the heart of understanding the problem
of conformal invariance. 

\begin{figure}
\centerline{\includegraphics[width=0.48\textwidth]{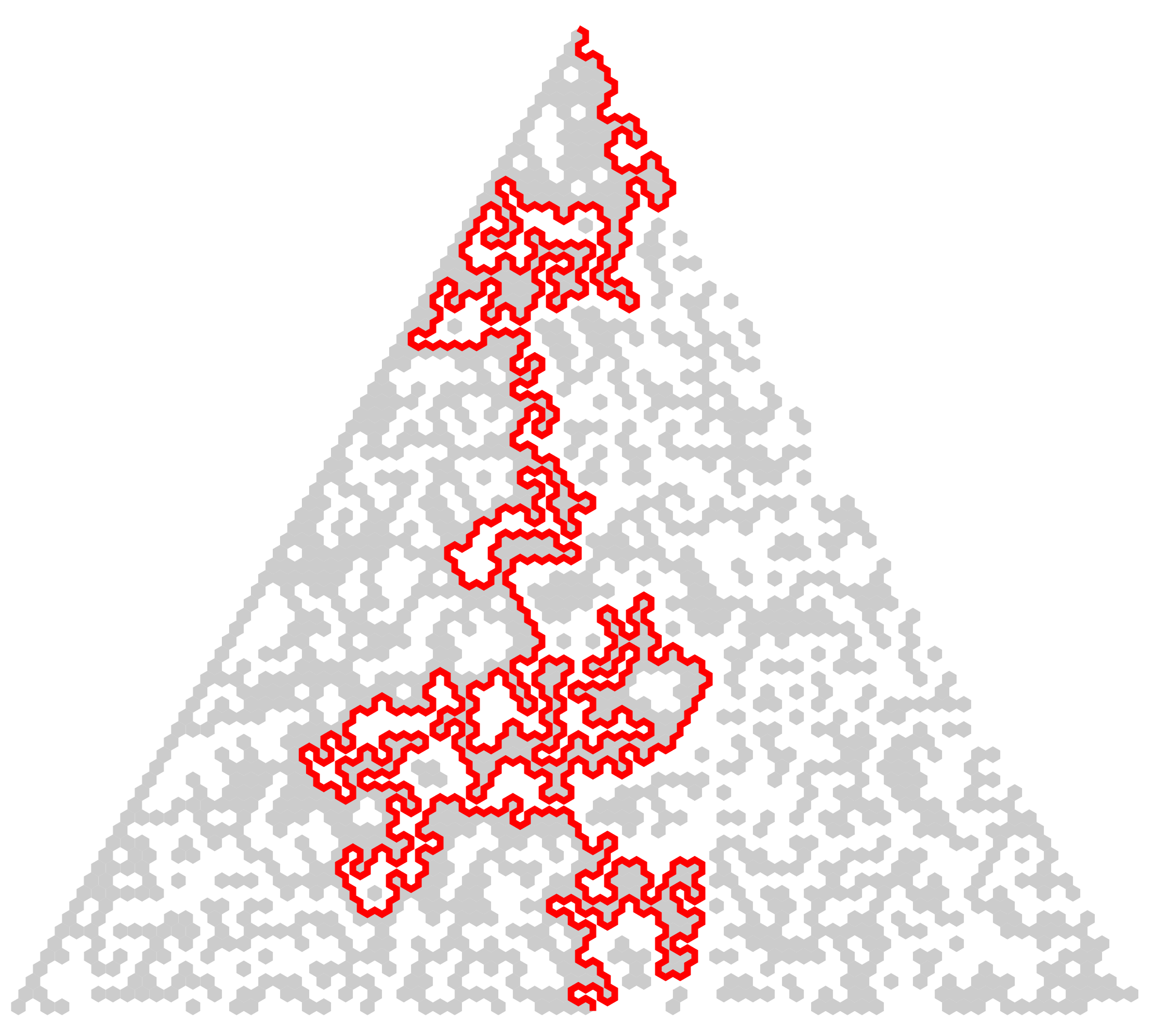}}
\caption*{{\bf Figure 8.} 
Site percolation on the triangular lattice $\TT$. The red path is a black/white interface.
Smirnov proved Cardy's formula, that 
states in this context that the hitting point of the interface on the bottom side is 
asymptotically uniformly distributed}
\end{figure}

Correlation inequalities are central to the theory of disordered systems in mathematics and physics.
The highly novel BK (van den Berg/Kesten) inequality plays a key role in systems subjected to a product measure
such as percolation. Proved in 1985, this inequality is a form of negative association,
based around the notion of the \lq disjoint occurrence'
of two events. It is a delicate and tantalizing result.

When $p>\pc$, there exists a.s.\ at least one infinite open cluster, but how many? This uniqueness problem
was answered by Harry in  joint  work with  Aizenman and  Newman. He tended to downplay his part in this work,
but his friends knew him better than to take such protestations at face value. Their paper
was soon superseded by the elegant argument of  Burton and  Keane, but it remains important
as a source of quantitative estimates.  

Under the title \lq ant in a labyrinth',  de Gennes proposed the use of a random walk to explore
the geometry of an open cluster. In a beautiful piece of work \cite{subdiffusive}, Harry showed the existence
of a measure known as the \lq  incipent infinite cluster',
obtained effectively by conditioning on critical percolation possessing an infinite cluster at the origin. 
He then proved that random walk $X_n$ on this cluster is subdiffusive in that there exists $\epsilon>0$ such that
$X_n/n^{\frac12-\epsilon} \to 0$.  This is in contrast to the situation on $\ZZ^d$ 
for which $X_n$ has order $n^{\frac12}$.   This \lq slowing down' occurs because the walk spends
time in blind alleys of the incipient infinite cluster. 

Whereas it was not possible to obtain exact results for critical percolation, Harry 
gave a  precise solution to the corresponding problem on the family tree $T$ of a critical branching process,
conditioned on non-extinction. Here also there are blind alleys, but it is possible to estimate the time spent in them. 
It turns out that
the displacement $X_n$ of the walk has order $n^{\frac13}$ and, moreover, 
Harry computed the limit distribution of $X_n/n^{\frac13}$.

In joint work with  Grimmett and Zhang on supercritical percolation, 
he showed that random walk on the infinite cluster in $d\ge 2$ dimensions  is recurrent if and only if $d=2$.
This basic result stimulated others to derive precise heat kernel estimates for random walks on percolation clusters and other 
random networks.  

We mention one further result for classical percolation on $\ZZ^d$. 
When $p<\pc$,  the tail of $|C|$ decays
exponentially to $0$, in that $\PP_p(|C| = n) \le e^{-\alpha n}$ for
some $\alpha(p)>0$. Matters are more complicated when $p>\pc$, since large clusters \lq prefer' to be infinite.
It turns out that the interior of a large cluster (of size $n$, say)
typically resembles that of the infinite cluster, and its
finiteness is controlled by its boundary (of order $n^{(d-1)/d}$).  
As a result, $|C|$ should have a tail of order $\exp\{-cn^{(d-1)/d}\}$.
There was a  proof by Aizenman,  Delyon, and  Souillard
that 
$$
\PP_p(|C|=n) \ge \exp\{-\beta n^{(d-1)/d}\}
$$ 
for $\beta(p)>0$.
Kesten and Zhang showed, by a block argument, 
the complementary inequality 
$$
\PP_p(|C|=n) \le  \exp\{-\gamma n^{(d-1)/d}\}
$$ 
for 
some $\gamma(p)>0$. 
When $d\ge 3$, they were in fact only able to show this for $p$ exceeding a certain value $\pslab$, but 
the conclusion for $p>\pc$ followed once  Grimmett and  Marstrand had proved 
the slab limit  $\pc=\pslab$.
Sharp asymptotics were established later by  Alexander,  Chayes, and  Chayes
when $d=2$, and by  Cerf
in the more challenging  situation of $d = 3$, in their work on the Wulff construction for percolation. 

\begin{figure}
\centerline{\includegraphics[width=0.45\textwidth]{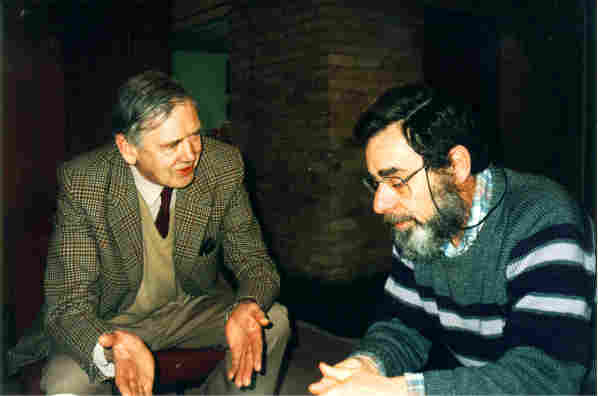}}
\caption*{{\bf Figure 9.}
Harry Kesten with John Hammersley, Oxford, 1993.}
\end{figure} 

First-passage percolation was introduced by Hammersley and Welsh 
in 1965  as an extension of classical percolation
in which each edge has a random \lq passage time', and one studies the set of vertices reached
from the origin along paths of  length not exceeding a given value. This is where
the notion of a subadditive stochastic process was introduced, and an ergodic theorem first proved.
The theory of subadditivity was useful throughout Harry's work on disordered networks, and 
indeed he noted that it provided an \lq\lq elegant''  proof of his 1960 theorem with Furstenberg on products of random matrices. 

Harry turned 
towards first-passage percolation around 1979, and he resolved a number of open problems,
and posed others, in a series of papers spanning nearly 10 years.  He established  fundamental properties
of the time constant, including positivity under a natural condition, and continuity as 
a function of the underlying distribution (with  Cox), together with a large deviation theorem for 
passage times (with Grimmett). Perhaps his most notable contribution was a theory of duality in three dimensions
akin to Whitney duality of two dimensions, as expounded in his St Flour notes, \cite{MR2905932}. 
The dual process is upon plaquettes, and dual \emph{surfaces} occupy the role of 
dual \emph{paths} in two dimensions. This leads to some tricky geometrical issues concerning  the combinatorics and topology
of dual surfaces which have been largely answered since by  Zhang,  Rossignol, 
Th\'eret,  Cerf, and others.

Percolation and its cousins provide the environment for a number of related processes
to which Harry contributed substantial results.
He was an enthusiastic contributor to too many collaborative ventures to be described in full  here.
As a sample we mention word percolation (with Benjamini,  Sidoravicius, Zhang),
$\rho$-percolation (with Su), random lattice animals (with Cox, Gandolfi, Griffin), 
and uniform spanning trees (with Benjamini, Peres, Schramm).

\section*{Postscript}

This memoir describes only a sample of Harry Kesten's impact 
within the probability and statistical physics communities.  Those seeking more may read  
Rick Durrett's account \cite{MR1703122}. The more comprehensive article \cite{grimmett}
includes summaries of Harry's work in other areas, such as
quasi-stationary distributions of Markov chains and bounded remainder sets in Euclidean dynamics.

The authors are grateful to several people for their advice:
Michael Kesten for access to the Kesten family papers
and for his personal reminiscences, Laurent Saloff-Coste for biographical material from Cornell,
Frank den Hollander for advice on Dutch matters, Rob van den Berg for 
his suggestions,
Hillel Furstenberg for his comments on a draft, and three anonymous readers for their suggestions. 

\begin{figure}[h]
\centerline{\includegraphics[width=0.205\textwidth]{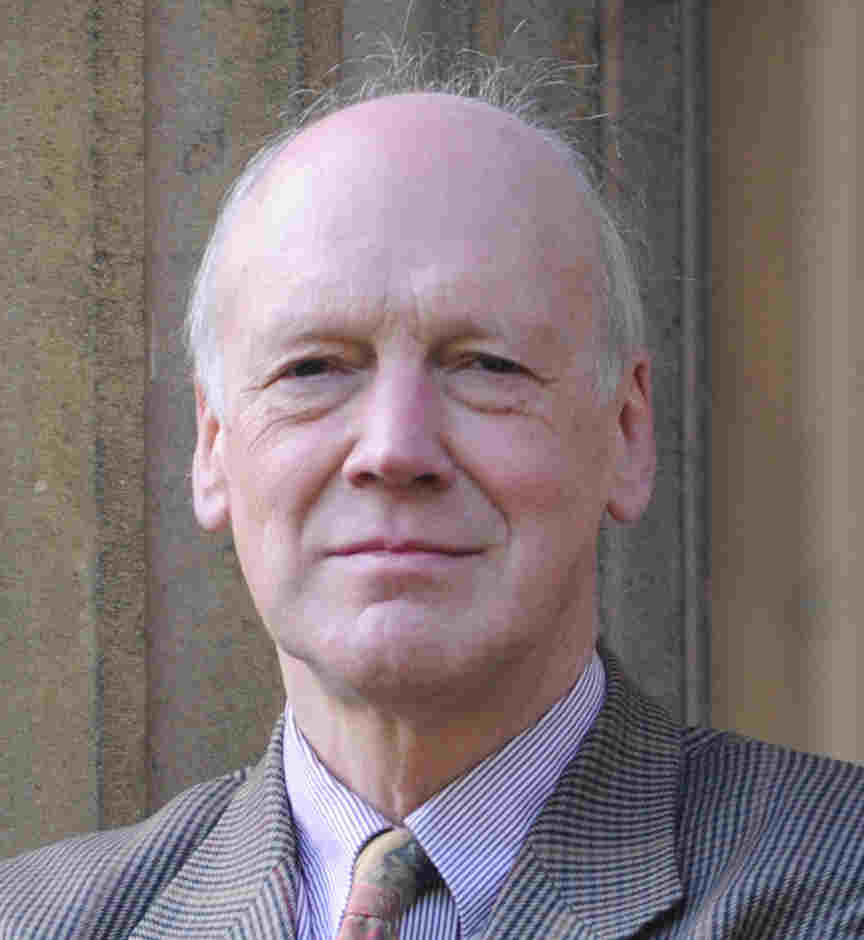}\quad
\includegraphics[width=0.22\textwidth]{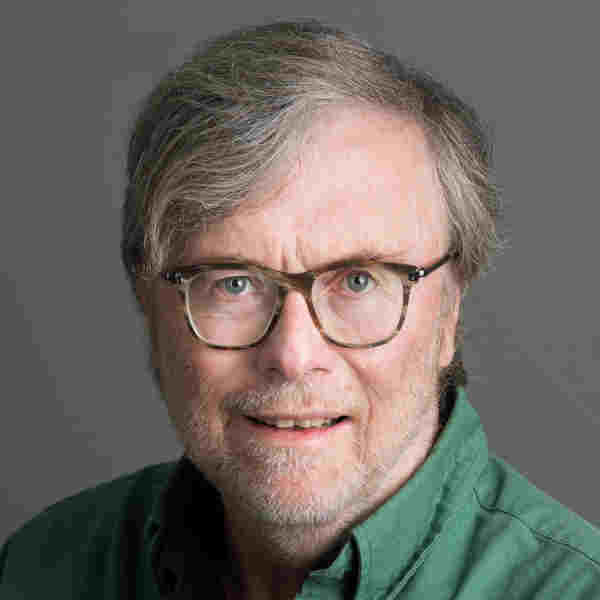}}
\caption*{{\bf Figure 10.}
The authors: Geoffrey Grimmett and Greg Lawler.}
\end{figure}


\subsection*{Credits}
Figure 1 is by permission of Cornell University.

\noindent
Figure 2 is by permission of Rob van den Berg.

\noindent
Figure 3 is by permission of Hugo Grimmett.

\noindent
Figure 4. Photographer: Konrad Jacobs. Source: Archives of the Mathematisches
Forschungsinstitut Oberwolfach.

\noindent
Figure 5 is by permission of Michael A.\ Morgan.

\noindent
Figures 6, 7, 8, 9 are by permission of Geoffrey Grimmett.

\noindent
Figure 10 is by permission of Downing College, Cambridge, and 
Jean Lachat/University of Chicago,
respectively.

\end{document}